# Algorithms for Construction, Classification and Enumeration of Closed Knight's Paths

*Stoyan Kapralov*
Technical University
of Gabrovo
Gabrovo, Bulgaria
s.kapralov@tugab.bg
0000-0002-5713-1488

*Valentin Bakoev*
"St. Cyril and St. Methodius"
University of Veliko Tarnovo,
Veliko Tarnovo, Bulgaria
v.bakoev@ts.uni-vt.bg
0000-0003-2503-5325

*Kaloyan Kapralov*
Sofia, Bulgaria
kaloyan.kapralov@gmail.com
0000-0003-4971-9422

***Abstract.*** Two algorithms for construction of all closed knight's paths of lengths up to 16 are presented. An approach for classification (up to equivalence) of all such paths is considered. By applying the construction algorithms and classification approach, we enumerate both unrestricted and non-intersecting knight's paths, and show the obtained results.

1. **Introduction.**

In graph theory, a *knight grap*h is a graph that represents all legal moves of a **chess knight** on a chessboard. An $m \times n$ knight graph is that for an $m \times n$ chessboard. An *open* knight's path visits every square on the chessboard exactly once and so it is a *Hamiltonian path* in the corresponding knight graph. When the starting point coincides with the endpoint, the knight's path is called a *closed* knight's path. It is a *Hamiltonian cycle* in this knight graph. The types of knight graphs and their properties, such as, for instance, that an $n \times n$ knight graph contains $4(n-2)(n-1)$ edges are discussed in [9].

The search for Hamiltonian paths and cycles and solving the many related problems has a long history – see [7] and the excellent survey presented on the internet in [8]. In the late 17th and early 18th centuries, mathematicians such as De Montmort and De Moivre provided solutions to the closed knight's path problem on the classical 8×8 board. But the first scientist who proposed a method for solving this problem and published an article about it was the famous Leonard Euler [1]. After him, well-known mathematicians such as Legendre, Vandermonde, and others also provided solutions to this problem. Special attention deserves the ingenious solution to the knight's path problem proposed by Warnsdorff in 1823. To get a solution, he suggests following a very simple rule: "Starting from an arbitrary square, always move the knight to an unvisited square from which there are as few subsequent moves as possible." This rule is the basis of a greedy algorithm that runs in linear time $O(n^2)$. Algorithms based on the "divide and conquer" strategy having the same type of time complexity are considered in [6]. Such types of approaches and algorithms are also discussed, starting with Euler's ideas. Although the problem of finding a Hamiltonian cycle is NP-complete in the general case, the knight graph has a special structure and such efficient solutions are possible. However, the problem of finding all nonequivalent (with respect to rotations and symmetries) solutions is very hard. That is why it was solved successfully only at the end of the 20th century [5].

In this paper, we consider a slightly different problem: to find, classify and enumerate all Hamiltonian cycles (i.e. undirected closed knight's paths) of definite lengths, both unbounded and non-self-intersecting (when viewed as polygons). Two cycles are called *equivalent* if one can be obtained from another by applying operations (one or more) translation, rotation, and symmetry on the chessboard, otherwise, they are *nonequivalent*. Here we classify and enumerate all nonequivalent cycles of lengths 12, 14, and 16. This problem has been solved for smaller cycle length values (6, 8 and 10) in our previous publications [2, 3].

In Section 2 we outline two algorithms for the construction of closed knight's paths and experimentally compare their execution times. In Section 3 we describe, via a concrete example, the equivalences – what it means for two paths to be equivalent, and for a given solution to be minimal. New results on the number of non-self-intersecting cycles are shown in Section 4. In the last section, we discuss some methodological aspects of the problem under consideration and its use in teaching subjects such as "Algorithms and Data Structures", "Competitive Programming", etc.

## 2. Algorithms for Construction of Closed Knight's Paths

Each closed knight's path has an even length, say $k = 2l$, since the knight's steps alternate the colors of the cells it steps on, for example black, white, black, white, …, black, or conversely: white, black, white, black, ..., white (an odd number of vertices in a path means an even number of edges, i.e. even length). Hence such a path can be embedded in a board of size $(k + 1) \times (k + 1)$. So we take an $m \times n$ board with $m = n = k + 1$ and number the cells with integers from 1 to $(k + 1)^2$, row-wise. We apply two similar approaches as outlined below.

In Algorithm 1, for each $s$ from the set $S$ of starting cells, we run a modified version of a depth-first search. Thus, via exhaustive backtracking, we build all knights' paths of length $k$ (with $k$ vertices and $k - 1$ edges). Finally, we check whether the last square is adjacent (via a chess knight's move) to the starting one, and if it is, that means that a new solution is found. Such a solution is saved if it is lexicographically minimal – the description of minimality and minimality check is given in the next section. Here is the C++ code of Algorithm 1.

```
void dfs()
{
  int u = p[len-1]; // p – list of the path cells
  if(len < k)
  {
    for (int v:adj[u])
      if(!vis[v])
      {
        p[len++] = v;
        vis[v] = true;
        dfs();
        len--;
        vis[v] = false;
      }
```

```
      }
      else
        if(adjacent(u,s))
        {
          path pt(p,p+k);
          if(is_min(pt))
            sol.push_back(pt);
        }
}

int main()
{
...
  for(s:S)
  {
    p[0] = s;
    vis[s] = true;
    len = 1;
    dfs();
  }
...
}
```

In Algorithm 2, for each pair of starting cell $s$ and ending cell $t$, we construct all connecting paths of length $k/2$. We then try to assemble a whole path of length $k$ from two halves. Of course, we may concatenate two paths if they do not intersect, that is, they have no common point except the two endpoints. It is sufficient to consider only cycles, which cannot be translated up or left, meaning that they have at least one point on the top row of the board and at least one point in the leftmost column. Thus, we define the equivalence between two cycles with respect to translation, i.e., sliding on the board. Hence, as well as for symmetry reasons, we may assume that we have to consider only the set $S = \{1,2,\ldots,k/2+1\}$ as a set of possible starting points. A comparison between the speed of Algorithm 1 (slower) and Algorithm 2 (faster) is shown in Table 1. The comparison was made on a common laptop with an 8th-generation Intel core i5 processor. Programs are written in C++.

| $k$ | Number of nonequivalent cycles | Algorithm 1 | Algorithm 2 |
| --- | --- | --- | --- |
| 6 | 25 | < 1 s | < 1 s |
| 8 | 480 | < 1 s | < 1 s |
| 10 | 12000 | 2 s | 1 s |
| 12 | 350256 | 44 s | 18 s |
| 14 | 10780549 | 2097 s | 692 s |
| 16 | 344680960 | ??? | ≈ 24 hours |

**Table 1.** The number of nonequivalent closed knight's paths
and the elapsed time in seconds of the two algorithms

## 3. Equivalent solutions

Let's take a look at the figures in Table 2. They are all obtained from the main figure in the upper left corner by rotating (counter-clockwise) in steps of 90° and applying symmetry about the vertical axis. We say that these 8 figures are equivalent to each other. Our task to enumerate and classify the closed knight's paths of a given length is essentially to enumerate some set of figures up to equivalence, i.e. to give just one representative from each equivalence class.

For pruning the unnecessary solutions we apply the following approach, illustrated in Table 2 for paths of length 8. The fields on the board are numbered row by row with the numbers $1, 2, \ldots, (k+1)^2$. In our example, the cells are numbered from 1 to 25 row-wise. Then every figure is represented by a sequence of 8 cell numbers. Each individual geometric figure represents 16 different sequences (8 ways to choose the starting point and two ways to choose the direction of the figure – clockwise or counter-clockwise). It is easy to see that for each of the 8 figures in Table 2, the lexicographically smallest of these 16 sequences is chosen to represent the corresponding figure. The figure in the upper left corner is minimal because its representing sequence is lexicographically smaller or equal to the representing sequences of the other 7 figures. Each of the other 7 figures is not minimal, because with a finite number of rotations and symmetries, it can be transformed into a figure with a lexicographically smaller representing sequence.

In the process of generating the figures with Algorithm 1 or Algorithm 2, if a figure that is not minimal is obtained, it is not taken into account, because it is equivalent to the corresponding minimal figure.

| 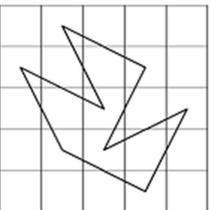 | 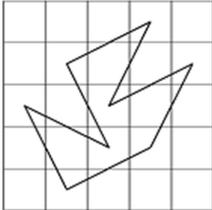 | 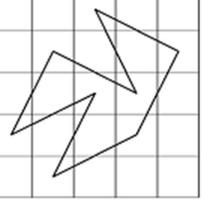 | 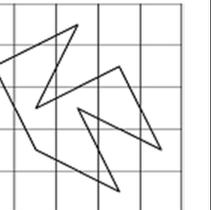 |
|---|---|---|---|
| 2,9,18,15,24,17,6,13 | 4,7,18,11,22,19,10,13 | 3,10,19,22,13,16,7,14 | 3,6,17,24,13,20,9 |
| 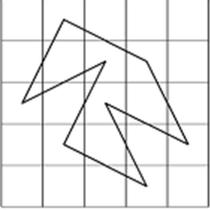 | 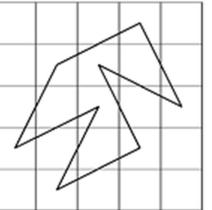 | 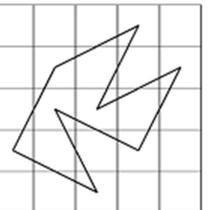 | 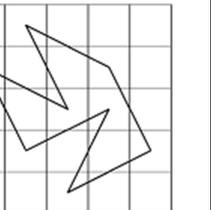 |
| 2,9,20,13,24,17,8,11 | 4,7,16,13,22,19,8,15 | 4,7,16,23,12,19,10,13 | 2,9,20,23,14,17,6,13 |

**Table 2.** The 8 equivalent solutions

In this way, there is no need to check if a newly constructed solution is equivalent to any of the already obtained thousands of solutions. It is sufficient just to check if it is a minimal one.

In [3] we began to distinguish closed knight's paths that contain the same sets of visited cells but the corresponding polygons are not congruent. We call such paths *geometrically distinct* closed knight's paths. The first such case occurs for length $k = 8$ and is illustrated in Table 3.

| 3 7 15 11 24 28 20 16 | 3 7 15 28 20 16 24 11 | 3 7 20 28 15 11 24 16 |

Table 3. Three geometrically distinct closed knight's paths of length 8

## 4. New results

An important problem related to the problem of all knight paths (open or closed) in an $m \times n$ chessboard is the problem of finding the uncrossed (i.e., non-self-intersecting) knight paths of maximum length [4, 8]. This problem was proposed by T. R. Dowson in 1930, who gave two solutions of length 35 for an $8 \times 8$ chessboard, and he later claimed (without proof) that they were of maximum length [4]. More about the solutions of the task in its various variants and the current results can be seen on the site [8].

With a suitable modification of the algorithms, we enumerated the non-equivalent and geometrically distinct closed knight's paths as well as all non-self-intersecting ones among them. Table 4 presents all currently known results for the number of closed knight's paths of length $k$, as well as the number of non-self-intersecting closed paths. New ones are marked with '*'. The previously known results are from [2] and [3].

| Length $k$ | Number of nonequivalent cycles | Number of non-self-intersecting cycles |
|---|---|---|
| 4 | 3 | 3 |
| 6 | 25 [2] | 13* |
| 8 | 480 [3] | 178* |
| 10 | 12000 [3] | 3034* |
| 12 | 350256* | 64877* |
| 14 | 10780549* | 1503790* |
| 16 | 344680960* | 36930111* |

Table 4. The new results

## 5. Conclusion

So far, the problem under consideration has been presented as a research problem. It was inspired by an assignment proposed by the first author at the 2015 National Spring Tournament in Informatics [10]. It is well known that the *knight's walk problem* and the *8 queens problem* are classic examples in the study of recursion and backtracking. That is why they occupy a well-deserved place in the relevant books or textbooks. So, the general problem of finding all open or closed knight's paths on an $m \times n$ chessboard and related problems have a pedagogical meaning. This general problem offers many possibilities and relevant examples in teaching subjects like "Programming", "Algorithms and data structures", "Graph algorithms", "Competitive programming", etc. It can be included in the study of topics like recursion, exhaustive search and backtracking, depth-first search, greedy algorithms, etc. For each of these subjects, the theme of the knight's paths or cycles can be used for lectures, labs, homework, and projects. Their difficulty can have different levels:

- initial – search for only one knight's path or cycle, experiments with chess boards of different sizes, square and rectangular;
- intermediate – search and enumerate all knight's paths/cycles, checking for the non-existence of a cycle in boards of certain sizes, search for semi-magic knight's paths/cycles (for details see [8]), etc.;
- high – learning and applying as a computer programs some equivalences and rejecting equivalent solutions which can be added to some intermediate-level problems;
- challenge – to check some results presented here and/or get new results.

Finally, we note that the knight's path problem continues to evolve, mostly due to the power of the latest processors. The site of G. Jelliss shows a lot of recent solutions to the various types of the knight's path problem in $m \times n$ chess boards, the size of which is increasing more and more [8].

## Acknowledgment

This work was supported in part by Grant No. 2209E/2022 of the Technical University of Gabrovo and by Veliko Tarnovo University Project FSD-31-243-23/21.03.2023.